\documentclass[11pt,a4paper]{article}

\usepackage[T1]{fontenc}
\usepackage{epsfig}
\usepackage{amsmath, amsthm}
\usepackage{amsfonts,amssymb}

\newtheorem{conj}{Conjecture}[section]
\newtheorem{theo}[conj]{Theorem}

\newtheorem{lem}[conj]{Lemma}
\newtheorem{prop}[conj]{Proposition}
\newtheorem{coro}[conj]{Corollary}

\newcommand{\de}{\mathrm{d}}
\newcommand{\e}{\mathrm{e}}

\newcommand{\Hess}{\mathrm{Hess}}

\newcommand{\R}{\mathbb{R}}
\newcommand{\N}{\mathbb{N}}

\newcommand{\eps}{\varepsilon}

\begin{document}

\title{A note on an $L^p$-Brunn-Minkowski inequality for convex measures in the unconditional case}

\author{Arnaud Marsiglietti}

\date{}

\maketitle
{\raggedright Institute for Mathematics and its Applications, University of Minnesota} \\

\noindent
arnaud.marsiglietti@ima.umn.edu \\

{\raggedright {\scriptsize This research was supported in part by the Institute for Mathematics and its Applications with funds provided by the National Science Foundation.}}

\begin{abstract}

We consider a different $L^p$-Minkowski combination of compact sets in $\R^n$ than the one introduced by Firey and we prove an $L^p$-Brunn-Minkowski inequality, $p \in [0,1]$, for a general class of measures called convex measures that includes log-concave measures, under unconditional assumptions. As a consequence, we derive concavity properties of the function $t \mapsto \mu(t^{\frac{1}{p}} A)$, $p \in (0,1]$, for unconditional convex measures $\mu$ and unconditional convex body $A$ in $\R^n$. We also prove that the (B)-conjecture for all uniform measures is equivalent to the (B)-conjecture for all log-concave measures, completing recent works by Saroglou.

\end{abstract}


\section{Introduction}

The Brunn-Minkowski inequality is a fundamental inequality in Mathematics, which states that for every convex subset $A, B \subset \R^n$ and for every $\lambda \in [0,1]$, one has
\begin{eqnarray}
|(1-\lambda)A + \lambda B|^{\frac{1}{n}} \geq (1-\lambda) |A|^{\frac{1}{n}} + \lambda |B|^{\frac{1}{n}},
\end{eqnarray}
where
$$ A+B = \{a + b; a \in A, b \in B \} $$
denotes the {\it Minkowski sum} of $A$ and $B$ and where $|\cdot|$ denotes the Lebesgue measure. The book by Schneider~\cite{S} and the survey by Gardner~\cite{G} famously reference the Brunn-Minkowski inequality and its consequences.

Several extensions of the Brunn-Minkowski inequality have been developed during the last decades by establishing functional versions (see $e.g.$ \cite{HM}, \cite{D}, \cite{DU}, \cite{U}), by considering different measures (see $e.g.$ \cite{B2}, \cite{B}), by generalizing the Minkowski sum (see $e.g.$ \cite{F}, \cite{F2}, \cite{F3}, \cite{Lu}, \cite{Lu2}), among others.  \\

In this paper, we will combine these extensions to prove an $L^p$-Brunn-Minkowski inequality for a large class of measures, including the log-concave measures.

Firstly, let us consider measures other than the Lebesgue measure. Following Borell~\cite{B2}, \cite{B}, we say that a Borel measure $\mu$ in $\R^n$ is {\it $s$-concave}, $s \in [-\infty, +\infty]$, if the inequality
$$ \mu((1 - \lambda)A + \lambda B) \geq M_s^{\lambda}(\mu(A),\mu(B)) $$
holds for every $\lambda \in [0,1]$ and for every compact subset $A,B \subset \mathbb{R}^n$ such that $\mu(A)\mu(B)>0$, where $M_s^{\lambda}(a,b)$ denotes the $s$-mean of the non-negative real numbers $a, b$ with weight $\lambda$, defined as
$$ M_s^{\lambda}(a,b) = ((1-\lambda)a^s + \lambda b^s)^{\frac{1}{s}} \quad \mbox{if $s \notin \{-\infty, 0, +\infty\}$}, $$
$M_{-\infty}^{\lambda}(a,b) = \min(a,b)$, $M_0^{\lambda}(a,b) = a^{1-\lambda} b^{\lambda}$, $M_{+\infty}^{\lambda}(a,b) = \max(a,b)$. Hence the Brunn-Minkowski inequality tells us that the Lebesgue measure in $\R^n$ is $\frac{1}{n}$-concave.

As a consequence of the H\"older inequality, one has $M_p^{\lambda}(a,b) \leq M_q^{\lambda}(a,b)$ for every $p \leq q$. Thus every $s$-concave measure is $- \infty$-concave. The $-\infty$-concave measures are also called {\it convex measures}. 

For $s \leq \frac{1}{n}$, Borell showed that every measure $\mu$, which is absolutely continuous with respect to the $n$-dimensional Lebesgue measure, is $s$-concave if and only if its density is an $\alpha$-concave function, with
\begin{eqnarray}\label{relation}
\alpha = \frac{s}{1-sn}  \in [-\frac{1}{n}, + \infty],
\end{eqnarray}
where a function $f : \R^n \to \R_+$ is said to be \textit{$\alpha$-concave}, with $\alpha \in [-\infty, +\infty]$, if the inequality
$$ f((1 - \lambda)x + \lambda y) \geq M_{\alpha}^{\lambda}(f(x), f(y)) $$
holds for every $x,y \in \R^n$ such that $f(x)f(y) > 0$ and for every $\lambda \in [0,1]$.

Secondly, let us consider a generalization of the notion of the Minkowski sum introduced by Firey, which leads to an {\it $L^p$-Brunn-Minkowski theory}. For {\it convex bodies} $A$ and $B$ in $\R^n$ ($i.e.$ compact convex sets containing the origin in the interior), the {\it $L^p$-Minkowski combination}, $p \in [-\infty, + \infty]$, of $A$ and $B$ with weight $\lambda \in [0,1]$ is defined by
$$ (1-\lambda) \cdot A \oplus_p \lambda \cdot B = \{x \in \R^n ; \langle x,u \rangle \leq M_p^{\lambda}(h_A(u), h_B(u)), \forall u \in S^{n-1} \}, $$
where $h_A$ denotes the {\it support function} of $A$ defined by
$$ h_A(u) = \max_{x \in A} \langle x,u \rangle, \quad u \in S^{n-1}. $$
Notice that for every $p \leq q$,
$$ (1-\lambda) \cdot A \oplus_p \lambda \cdot B \subset (1-\lambda) \cdot A \oplus_q \lambda \cdot B. $$
The support function is an important tool in Convex Geometry, having the property to determine a convex body and to be linear with respect to Minkowski sum and dilation:
$$ A = \{x \in \R^n ; \langle x,u \rangle \leq h_A(u), \forall u \in S^{n-1}\}, \quad h_{A+B} = h_A + h_B, \quad h_{\mu A} = \mu h_A, $$
for every convex body $A,B$ in $\R^n$ and every scalar $\mu \geq 0$. Thus,
$$ (1-\lambda) \cdot A \oplus_1 \lambda \cdot B = (1-\lambda)A + \lambda B. $$

In this paper, we consider a different $L^p$-Minkowski combination. Before giving the definition, let us recall that a function $f : \R^n \to \R$ is {\it unconditional} if there exists a basis $(a_1, \cdots, a_n)$ of $\R^n$ (the canonical basis in the sequel) such that for every $x=\sum_{i=1}^n x_i a_i \in \R^n$ and for every $\eps=(\eps_1, \cdots,\eps_n) \in \{-1,1\}^n$, one has $f(\sum_{i=1}^n \eps_i x_i a_i) = f(x)$. For $p=(p_1, \cdots, p_n) \in [-\infty, +\infty]^n$, $a=(a_1, \cdots, a_n) \in (\R_+)^n, b=(b_1, \cdots, b_n) \in (\R_+)^n$ and $\lambda \in [0,1]$, let us denote
$$ (1-\lambda) a +_p \lambda b = (M_{p_1}^{\lambda}(a_1, b_1), \cdots, M_{p_n}^{\lambda}(a_n, b_n)) \in (\R_+)^n. $$
For non-empty subsets $A, B \subset \R^n$, for $p \in [-\infty, +\infty]^n$ and for $\lambda \in [0,1]$, we define the $L^p$-Minkowski combination of $A$ and $B$ with weight $\lambda$, denoted by $(1-\lambda) \cdot A +_p \lambda \cdot B$, to be the unconditional subset ($i.e.$ the indicator function is unconditional) such that
$$ ((1-\lambda) \cdot A +_p \lambda \cdot B) \cap (\R_+)^n = \{(1-\lambda) a +_p \lambda b ~;~ a \in A \cap (\R_+)^n, b \in B \cap (\R_+)^n \}. $$
This definition is consistent with the well known fact that an unconditional set (or function) is entirely determined on the positive octant $(\R_+)^n$. Moreover, this $L^p$-Minkowski combination coincides with the classical Minkowski sum when $p=(1, \cdots, 1)$ and $A,B$ are unconditional convex subsets of $\R^n$ (see Proposition~\ref{recover} below).

Using an extension of the Brunn-Minkowski inequality discovered by Uhrin~\cite{U}, we prove the following result:

\begin{theo}\label{BMI}

Let $p=(p_1, \cdots, p_n) \in [0,1]^n$ and let $\alpha \in \R$ such that $\alpha \geq - \left( \sum_{i=1}^n p_i^{-1} \right)^{-1}$. Let $\mu$ be an unconditional measure in $\R^n$ that has an $\alpha$-concave density function with respect to the Lebesgue measure. Then, for every unconditional convex body $A,B$ in $\R^n$ and for every $\lambda \in [0,1]$,
\begin{eqnarray}
\mu((1-\lambda) \cdot A +_p \lambda \cdot B) \geq M_{\gamma}^{\lambda} (\mu(A), \mu(B)),
\end{eqnarray}
where $\gamma = \left( \sum_{i=1}^{n} p_i^{-1} + \alpha^{-1} \right)^{-1}$.

\end{theo}

The case of the Lebesgue measure and $p=(0, \cdots, 0)$ is treated by Saroglou~\cite{S}, answering a conjecture by B\"or\"oczky, Lutwak, Yang and Zhang~\cite{BLYZ} in the unconditional case.

\begin{conj}[log-Brunn-Minkowski inequality \cite{BLYZ}]\label{log-BMI}

Let $A,B$ be symmetric convex bodies in $\R^n$ and let $\lambda \in [0,1]$. Then,
\begin{eqnarray}\label{log-BM}
|(1-\lambda) \cdot A \oplus_0 \lambda \cdot B| \geq |A|^{1-\lambda} |B|^{\lambda}.
\end{eqnarray}

\end{conj}

Useful links have been discovered by Saroglou~\cite{S}, \cite{S2} between Conjecture~\ref{log-BMI} and the (B)-conjecture.

\begin{conj}[(B)-conjecture \cite{L}, \cite{CFM}]\label{(B)}

Let $\mu$ be a symmetric log-concave measure in $\R^n$ and let $A$ be a symmetric convex subset of $\R^n$. Then the function $t \mapsto \mu(\e^t A)$ is log-concave on $\R$.

\end{conj}

The (B)-conjecture was solved by Cordero-Erausquin, Fradelizi and Maurey~\cite{CFM} for the Gaussian measure and for the unconditional case. As a variant of the (B)-conjecture, one may study concavity properties of the function $t \mapsto \mu(V(t) A)$ where $V : \R \to \R_+$ is a convex function. As a consequence of Theorem~\ref{BMI}, we deduce concavity properties of the function $t \mapsto \mu(t^{\frac{1}{p}} A)$, $p \in (0,1]$, for every unconditional $s$-concave measure $\mu$ and every unconditional convex body $A$ in $\R^n$ (see Proposition~\ref{(B)-variant} below).

Saroglou~\cite{S2} also proved that the log-Brunn-Minkowski inequality for the Lebesgue measure (inequality~(\ref{log-BM})) is equivalent to the log-Brunn-Minkowski inequality for all log-concave measures. We continue these kinds of equivalences by proving that the (B)-conjecture for all uniform measures is equivalent to the (B)-conjecture for all log-concave measures (see Proposition~\ref{reduction} below).

We also investigate functional versions of the (B)-conjecture, which may be read as follows:

\begin{conj}[Functional version of the (B)-conjecture]\label{function}

Let $f,g : \R^n \to \R_+$ be even log-concave functions. Then the function
$$ t \mapsto \int_{\R^n} f(\e^{-t}x) g(x) \, \de x $$
is log-concave on $\R$.

\end{conj}

We prove that Conjecture~\ref{function} is equivalent to Conjecture~\ref{(B)} (see Proposition~\ref{equiv} below). \\

Let us note that other developments in the use of the earlier mentioned extensions of the Brunn-Minkowski inequality have been recently made as well (see $e.g.$~\cite{Bo}, \cite{Ca}, \cite{C}, \cite{GHW}). \\

The rest of the paper is organized as follows: In the next section, we prove Theorem~\ref{BMI} and we extend it to $m$ sets, $m \geq 2$. We also compare our $L^p$-Minkowski combination to the Firey combination and derive an $L^p$-Brunn-Minkowski inequality for the Firey combination. We then discuss the consequences of a variant of the (B)-conjecture, namely we deduce concavity properties of the function $t \mapsto \mu(t^{\frac{1}{p}} A)$, $p \in (0,1]$. In Section~3, we prove that the (B)-conjecture for all uniform measures is equivalent to the (B)-conjecture for all log-concave measures, and we also prove that the (B)-conjecture is equivalent to its functional version Conjecture~\ref{function}.

\section{Proof of Theorem~\ref{BMI} and consequences}

\subsection{Proof of Theorem~\ref{BMI}}

Before proving Theorem~\ref{BMI}, let us show that our $L^p$-Minkowski combination coincides with the classical Minkowski sum when $p=(1, \cdots, 1)$, for unconditional convex sets.

\begin{prop}\label{recover}

Let $A,B$ be unconditional convex subsets of $\R^n$ and let $\lambda \in [0,1]$. Then,
$$ (1-\lambda) \cdot A +_{\widetilde{1}} \lambda \cdot B = (1-\lambda)A + \lambda B, $$
where ${\widetilde{1}} = (1, \cdots, 1)$.

\end{prop}

\begin{proof}
Since the sets $(1-\lambda) \cdot A +_{\widetilde{1}} \lambda \cdot B$ and $(1-\lambda)A + \lambda B$ are unconditional, it is sufficient to prove that
$$ ((1-\lambda) \cdot A +_{\widetilde{1}} \lambda \cdot B) \cap (\R_+)^n = ((1-\lambda)A + \lambda B) \cap (\R_+)^n. $$

Let $x \in ((1-\lambda) A + \lambda B) \cap (\R_+)^n$. There exists $a=(a_1, \cdots, a_n) \in A$ and $b=(b_1, \cdots, b_n) \in B$ such that $x=(1-\lambda)a + \lambda b$ and for every $i \in \{1, \cdots, n\}$, $(1-\lambda)a_i + \lambda b_i \in \R_+$. Let $\eps, \eta \in \{-1, 1\}^n$ such that $(\eps_1 a_1, \cdots, \eps_n a_n) \in (\R_+)^n$ and $(\eta_1 b_1, \cdots, \eta_n b_n) \in (\R_+)^n$. Notice that for every $i \in \{1, \cdots, n\}$, $0 \leq (1-\lambda)a_i + \lambda b_i \leq (1-\lambda)\eps_i a_i + \lambda \eta_i b_i$. Since the sets $A$ and $B$ are convex and unconditional, it follows that $x \in (1-\lambda) (A \cap (\R_+)^n) + \lambda (B \cap (\R_+)^n) = ((1-\lambda) \cdot A +_{\widetilde{1}} \lambda \cdot B) \cap (\R_+)^n$.

The other inclusion is clear due to the definition of the set $(1-\lambda) \cdot A +_{\widetilde{1}} \lambda \cdot B$.
\end{proof}

\begin{proof}[Proof of Theorem~\ref{BMI}]
Let $\lambda \in [0,1]$ and let $A,B$ be unconditional convex bodies in $\R^n$. 

It has been shown by Uhrin~\cite{U} that if $f,g,h : (\R_+)^n \to \R_+$ are bounded measurable functions such that for every $x,y \in (\R_+)^n$, $ h((1-\lambda)x +_p \lambda y) \geq M_{\alpha}^{\lambda}(f(x), g(y))$, then
$$ \int_{(\R_+)^n} h(x) \, \de x \geq M_{\gamma}^{\lambda} \left( \int_{(\R_+)^n} f(x) \, \de x, \int_{(\R_+)^n} g(x) \, \de x \right), $$
where $\gamma = \left( \sum_{i=1}^n p_i^{-1} + \alpha^{-1} \right)^{-1}$.

Let us denote by $\phi$ the density function of $\mu$ and let us set $h=1_{(1-\lambda) \cdot A +_p \lambda \cdot B} \phi$, $f=1_A \phi$ and $g=1_B \phi$. By assumption, the function $\phi$ is unconditional and $\alpha$-concave, hence $\phi$ is non-increasing in each coordinate on the octant $(\R_+)^n$. Then for every $x,y \in (\R_+)^n$, one has
$$ \phi((1-\lambda)x +_p \lambda y) \geq \phi((1-\lambda)x + \lambda y) \geq M_{\alpha}^{\lambda}(\phi(x), \phi(y)). $$
Hence,
$$ h((1-\lambda)x +_p \lambda y) \geq M_{\alpha}^{\lambda}(f(x), g(y)). $$
Thus we may apply the result mentioned at the beginning of the proof to obtain that
$$ \int_{(\R_+)^n} h(x) \, \de x \geq M_{\gamma}^{\lambda} \left( \int_{(\R_+)^n} f(x) \, \de x, \int_{(\R_+)^n} g(x) \, \de x \right), $$
where $\gamma = \left( \sum_{i=1}^n p_i^{-1} + \alpha^{-1} \right)^{-1}$. In other words, one has
$$ \mu(((1-\lambda) \cdot A +_p \lambda \cdot B) \cap (\R_+)^n) \geq M_{\gamma}^{\lambda}(\mu(A \cap (\R_+)^n), \mu(B \cap (\R_+)^n)). $$
Since the sets $(1-\lambda) \cdot A +_p \lambda \cdot B$, $A$ and $B$ are unconditional, it follows that
$$ \mu((1-\lambda) \cdot A +_p \lambda \cdot B) \geq M_{\gamma}^{\lambda}(\mu(A), \mu(B)). $$
\end{proof}

\paragraph*{Remark.}
One may similarly define the $L^p$-Minkowski combination
$$ \lambda_1 \cdot A_1 +_p \cdots +_p \lambda_m \cdot A_m, $$
for $m$ convex bodies $A_1, \dots, A_m \subset \R^n$, $m \geq 2$, where $\lambda_1, \dots, \lambda_m \in [0,1]$ are such that $\sum_{i=1}^m \lambda_i = 1$, by extending the definition of the $p$-mean $M_p^{\lambda}$ to $m$ non-negative numbers. By induction, one has under the same assumptions of Theorem~\ref{BMI} that
\begin{eqnarray}\label{induction}
\mu(\lambda_1 \cdot A_1 +_p \cdots +_p \lambda_m \cdot A_m) \geq M_{\gamma}^{\lambda}(\mu(A_1), \cdots, \mu(A_m)),
\end{eqnarray}
where $\gamma = \left( \sum_{i=1}^{n} p_i^{-1} + \alpha^{-1} \right)^{-1}$. Indeed, let $m \geq 2$ and let us assume that inequality~(\ref{induction}) holds. Notice that
$$ \lambda_1 \cdot A_1 +_p \cdots +_p \lambda_m \cdot A_m +_p \lambda_{m+1} \cdot A_{m+1} = \left( \sum_{i=1}^{m} \lambda_i \right) \cdot \widetilde{A}  +_p \lambda_{m+1} \cdot A_{m+1}, $$
where
$$ \widetilde{A} := \left( \frac{\lambda_1}{\sum_{i=1}^{m} \lambda_i} \cdot A_1 +_p \cdots +_p \frac{\lambda_m}{\sum_{i=1}^{m} \lambda_i} \cdot A_m \right). $$
Thus,
\begin{eqnarray*}
\mu \left( \left( \sum_{i=1}^{m} \lambda_i \right) \cdot \widetilde{A}  +_p \lambda_{m+1} \cdot A_{m+1} \right) & \geq & \left( \left( \sum_{i=1}^{m} \lambda_i \right) \mu(\widetilde{A})^{\gamma} + \lambda_{m+1} \mu(A_{m+1})^{\gamma} \right)^{\frac{1}{\gamma}} \\ & \geq & \left( \sum_{i=1}^{m+1} \lambda_i \mu(A_i)^{\gamma} \right)^{\frac{1}{\gamma}}.
\end{eqnarray*}

\subsection{Consequences}

The following result compares the $L^p$-Minkowski combinations $\oplus_p$ and $+_p$.

\begin{lem}\label{include}

Let $p \in [0,1]$ and set $\widetilde{p}=(p, \cdots, p) \in [0,1]^n$. For every unconditional convex body $A,B$ in $\R^n$ and for every $\lambda \in [0,1]$, one has
$$ (1-\lambda) \cdot A \oplus_p \lambda \cdot B \supset (1-\lambda) \cdot A +_{\widetilde{p}} \lambda \cdot B. $$

\end{lem}

\begin{proof}
The case $p=0$ is proved in~\cite{S}. Let $p \neq 0$. Since the sets $(1-\lambda) \cdot A \oplus_p \lambda \cdot B$ and $(1-\lambda) \cdot A +_{\widetilde{p}} \lambda \cdot B$ are unconditional, it is sufficient to prove that
$$ ((1-\lambda) \cdot A \oplus_p \lambda \cdot B) \cap (\R_+)^n \supset ((1-\lambda) \cdot A +_{\widetilde{p}} \lambda \cdot B) \cap (\R_+)^{n}. $$
Let $u \in S^{n-1} \cap (\R_+)^n$ and let $x \in ((1-\lambda) \cdot A +_{\widetilde{p}} \lambda \cdot B) \cap (\R_+)^{n} $. One has,
\begin{eqnarray*}
\langle x,u \rangle & = & \sum_{i=1}^n ((1-\lambda)a_i^p + \lambda b_i^p)^{\frac{1}{p}} u_i \\ & = & \sum_{i=1}^n ((1-\lambda)(a_i u_i)^p + \lambda (b_i u_i)^p)^{\frac{1}{p}} \\ & = & \|(1-\lambda) X + \lambda Y\|_{\frac{1}{p}}^{\frac{1}{p}},
\end{eqnarray*}
where $X = ((a_1 u_1)^p, \cdots, (a_n u_n)^p)$ and $Y = ((b_1 u_1)^p, \cdots, (b_n u_n)^p)$. Notice that $\|X\|_{\frac{1}{p}} \leq h_A(u)^p$, $\|Y\|_{\frac{1}{p}} \leq h_B(u)^p$ and that $\|\cdot\|_{\frac{1}{p}}$ is a norm. It follows that
$$ \langle x,u \rangle \leq \left( (1-\lambda) \|X\|_{\frac{1}{p}} + \lambda \|Y\|_{\frac{1}{p}} \right)^{\frac{1}{p}} \leq \left( (1-\lambda) h_A(u)^p + \lambda h_B(u)^p \right)^{\frac{1}{p}}. $$
Hence, $x \in ((1-\lambda) \cdot A \oplus_p \lambda \cdot B) \cap (\R_+)^{n}$.
\end{proof}

From Lemma~\ref{include} and Theorem~\ref{BMI}, one obtains the following result:

\begin{coro}\label{Lp}

Let $p \in [0,1]$. Let $\mu$ be an unconditional measure in $\R^n$ that has an $\alpha$-concave density function, with $\alpha \geq - \frac{p}{n}$. Then for every unconditional convex body $A,B$ in $\R^n$ and for every $\lambda \in [0,1]$,
\begin{eqnarray}\label{LpI}
\mu((1-\lambda) \cdot A \oplus_p \lambda \cdot B) \geq M_{\gamma}^{\lambda} (\mu(A), \mu(B)),
\end{eqnarray}
where $\gamma = \left( \frac{n}{p} + \frac{1}{\alpha} \right)^{-1}$.

\end{coro}

\noindent
{\bf Remarks.}
\begin{enumerate}
\item By taking $\alpha = 0$ in Corollary~\ref{Lp} (corresponding to log-concave measures), one obtains
$$ \mu((1-\lambda) \cdot A \oplus_0 \lambda \cdot B) \geq \mu(A)^{1-\lambda} \mu(B)^{\lambda}. $$
\item By taking $\alpha = + \infty$ in Corollary~\ref{Lp} (corresponding to $\frac{1}{n}$-concave measures), one obtains that for every $p \in [0,1]$,
$$ \mu((1-\lambda) \cdot A \oplus_p \lambda \cdot B)^{\frac{p}{n}} \geq (1-\lambda) \mu(A)^{\frac{p}{n}} + \lambda \mu(B)^{\frac{p}{n}}.  $$
Equivalently, for every $p \in [0,1]$, for every unconditional convex body $A,B$ in $\R^n$ and for every unconditional convex set $K \subset \R^n$,
$$ |((1-\lambda) \cdot A \oplus_p \lambda \cdot B) \cap K|^{\frac{p}{n}} \geq (1-\lambda) |A \cap K|^{\frac{p}{n}} + \lambda |B \cap K|^{\frac{p}{n}}. $$
\end{enumerate}

Let us recall that the function $t \mapsto \mu(\e^t A)$ is log-concave on $\R$ for every unconditional log-concave measure $\mu$ and every unconditional convex body $A$ in $\R^n$ (see~\cite{CFM}). By adapting the argument of~\cite{M}, Proof of Proposition~3.1 (see Proof of Corollary~\ref{concave} below), it follows that the function $t \mapsto \mu(t^{\frac{1}{p}} A)$ is $\frac{p}{n}$-concave on $\R_+$, for every $p \in (0,1]$, for every unconditional $s$-concave measure $\mu$, with $s \geq 0$, and for every unconditional convex body $A$ in $\R^n$. However, no concavity properties are known for the function $t \mapsto \mu(\e^t A)$ when $\mu$ is an $s$-concave measure with $s < 0$. Instead, for these measures we prove concavity properties of the function $t \mapsto \mu(t^{\frac{1}{p}} A)$.

\begin{prop}\label{(B)-variant}

Let $p \in (0,1]$, let $\mu$ be an unconditional measure that has an $\alpha$-concave density function, with $\alpha \in [-\frac{p}{n}, 0)$ and let $A$ be an unconditional convex body in $\R^n$. Then the function $t \mapsto \mu(t^{\frac{1}{p}} A)$ is $\left( \frac{n}{p} + \frac{1}{\alpha} \right)^{-1}$-concave on $\R_+$.

\end{prop}

\begin{proof}
Let $t_1, t_2 \in \R_+$. By applying Corollary~\ref{Lp} to the sets $t_1^{\frac{1}{p}} A$ and $t_2^{\frac{1}{p}} A$, one obtains
$$ \mu(((1-\lambda) t_1 + \lambda t_2)^{\frac{1}{p}} A) = \mu((1-\lambda) \cdot t_1^{\frac{1}{p}} A \oplus_p \lambda \cdot t_2^{\frac{1}{p}} A) \geq M_{\gamma}^{\lambda} (\mu(t_1^{\frac{1}{p}} A), \mu(t_2^{\frac{1}{p}} A)), $$
where $\gamma = \left( \frac{n}{p} + \frac{1}{\alpha} \right)^{-1}$. Hence the function $t \mapsto \mu(t^{\frac{1}{p}} A)$ is $\gamma$-concave on $\R_+$.
\end{proof}

As a consequence, we derive concavity properties for the function $t \mapsto \mu(t A)$.

\begin{coro}\label{concave}

Let $p \in (0,1]$, let $\mu$ be an unconditional measure that has an $\alpha$-concave density function, with $\alpha \in [-\frac{p}{n}, 0)$, and let $A$ be an unconditional convex body in $\R^n$. Then the function $t \mapsto \mu(t A)$ is $\left(\frac{1-p}{n} + \gamma \right)$-concave on $\R_+$, where $\gamma = \left( \frac{n}{p} + \frac{1}{\alpha} \right)^{-1}$.

\end{coro}

\begin{proof}
We adapt~\cite{M}, Proof of Proposition~3.1. Let us denote by $\phi$ the density function of the measure $\mu$ and let us denote by $F$ the function $t \mapsto \mu(tA)$. From Proposition~\ref{(B)-variant}, the function $t \mapsto F(t^{\frac{1}{p}})$ is $\gamma$-concave, hence the right derivative of $F$, denoted by $F_+'$, exists everywhere and the function $t \mapsto \frac{1}{p} t^{\frac{1}{p} - 1} F_+'(t^{\frac{1}{p}}) F(t^{\frac{1}{p}})^{\gamma - 1}$ is non-increasing. Notice that
$$ F(t) = t^n \int_{A} \phi(tx) \, \de x, $$
and that $t \mapsto \phi(tx)$ is non-increasing, thus the function $t \mapsto \frac{1}{{t^{1-p}}} F(t)^{\frac{1-p}{n}}$ is non-increasing. Since
$$ F_+'(t)F(t)^{\frac{1-p}{n} + \gamma - 1} = t^{1-p} F_+'(t) F(t)^{\gamma - 1} \cdot \frac{1}{{t^{1-p}}} F(t)^{\frac{1-p}{n}}, $$
it follows that $F_+'(t)F(t)^{\frac{1-p}{n} + \gamma - 1}$ is non-increasing as the product of two non-negative non-increasing functions. Hence $F$ is $\left(\frac{1-p}{n} + \gamma \right)$-concave.
\end{proof}

\paragraph*{Remark.}
For every $s$-concave measure $\mu$ and every convex subset $A \subset \R^n$, the function $t \mapsto \mu(t A)$ is $s$-concave. Hence Corollary~\ref{concave} is of value only if $\frac{1-p}{n} + \gamma \geq \frac{\alpha}{1+\alpha n}$ (see relation~(\ref{relation})). Notice that this condition is satisfied if $\alpha \geq -\frac{p}{n(1+p)}$. We thus obtain the following corollary:

\begin{coro}\label{improve}

Let $p \in (0,1]$, let $\mu$ be an unconditional measure that has an $\alpha$-concave density function, with $-\frac{p}{n(1+p)} \leq \alpha < 0$ and let $K$ be an unconditional convex body in $\R^n$. Then, for every subsets $A, B \in \{\mu K; \mu > 0 \}$ and every $\lambda \in [0,1]$, one has
$$ \mu((1-\lambda)A + \lambda B) \geq M_{\frac{1-p}{n} + \gamma}^{\lambda}(\mu(A), \mu(B)), $$
where $\gamma = \left( \frac{n}{p} + \frac{1}{\alpha} \right)^{-1}$.

\end{coro}

In~\cite{M}, the author investigated improvements of concavity properties of convex measures under additional assumptions, such as symmetries. Notice that Corollary~\ref{improve} follows the same path and completes the results that can be found in~\cite{M}. Let us conclude this section by the following remark, which concerns the question of the improvement of concavity properties of convex measures.

\paragraph*{Remark.}

Let $\mu$ be a Borel measure that has a density function with respect to the Lebesgue measure in $\R^n$. One may write the density function of $\mu$ in the form $\e^{-V}$, where $V : \R^n \to \R$ is a measurable function. Let us assume that $V$ is $C^2$. Let $\gamma \in \R \setminus \{0\}$. The function $\e^{-V}$ is $\gamma$-concave if $\Hess(\gamma \e^{-\gamma V})$, the Hessian of $\gamma \e^{-\gamma V}$, is non-positive (in the sense of symmetric matrices). One has
$$ \Hess(\gamma \e^{-\gamma V}) = - \gamma^2 \nabla \cdot (\nabla V \e^{-\gamma V}) = \gamma^2 \e^{-V} (\gamma \nabla V \otimes \nabla V -\Hess(V)), $$
where  $\nabla V \otimes \nabla V = \left( \frac{\partial V}{\partial x_i} \frac{\partial V}{\partial x_j} \right)_{1 \leq i,j \leq n}$. Hence the matrix $\Hess(\gamma \e^{-\gamma V})$ is non-positive if and only if the matrix $\gamma \nabla V \otimes \nabla V - \Hess(V)$ is non-positive. \\

Let us apply this remark to the Gaussian measure
$$ \de \gamma_n(x) = \frac{1}{(2 \pi)^{\frac{n}{2}}} \e^{-\frac{|x|^2}{2}} \, \de x, \quad x \in \R^n. $$
Here $V(x) = \frac{|x|^2}{2} + c_n$, where $c_n = \frac{n}{2} \log(2 \pi)$. Thus, $\nabla V \otimes \nabla V = (x_i x_j)_{1 \leq i,j \leq n}$ and $\Hess(V) = Id$ the Identity matrix. Notice that the eigenvalues of $\gamma \nabla V \otimes \nabla V - \Hess(V)$ are $-1$ (with multiplicity $(n-1)$) and $\gamma |x|^2 - 1$. Hence if $\gamma |x|^2 - 1 \leq 0$, then $\gamma \nabla V \otimes \nabla V - \Hess(V)$ is non-positive. One deduces that for every $\gamma > 0$, for every compact sets $A,B \subset \frac{1}{\sqrt{\gamma}} B_2^n$ and for every $\lambda \in [0,1]$, one has
\begin{eqnarray}\label{improvement}
\gamma_n((1-\lambda)A + \lambda B) \geq M_{\frac{\gamma}{1 + \gamma n}}^{\lambda}(\gamma_n(A), \gamma_n(B)),
\end{eqnarray}
where $B_2^n$ denotes the Euclidean closed unit ball in $\R^n$.

Since the Gaussian measure is a log-concave measure, inequality~(\ref{improvement}) is an improvement of the concavity of the Gaussian measure when restricted to compact sets $A,B \subset \frac{1}{\sqrt{\gamma}} B_2^n$.

\section{Equivalence between (B)-conjecture-type problems}

In the following proposition, we demonstrate that it is sufficient to prove the (B)-conjecture for all uniform measures in $\R^n$, for every $n \in \N^*$, to obtain the (B)-conjecture for all symmetric log-concave measures in $\R^n$, for every $n \in \N^*$. This completes recent works by Saroglou~\cite{S}, \cite{S2}.

In the following, we say that a measure $\mu$ satisfies the (B)-property if the function $t \mapsto \mu(\e^t A)$ is log-concave on $\R$ for every symmetric convex set $A \subset \R^n$.

\begin{prop}\label{reduction}

If every symmetric uniform measure in $\R^n$, for every $n \in \N^*$, satisfies the (B)-property, then every symmetric log-concave measure in $\R^n$, for every $n \in \N^*$, satisfies the (B)-property.

\end{prop}

\begin{proof}
The proof is inspired by~\cite{AKM}, beginning of Section 3. \\

{\bf Step 1:} \underline{Stability under orthogonal projection} \\
Let us show that the (B)-property is stable under orthogonal projection onto an arbitrary subspace. 

Let $F$ be a $k$-dimensional subspace of $\R^n$. Let us define for every compactly supported measure $\mu$ in $\R^n$ and every measurable subset $A \subset F$,
$$ \Pi_F \mu (A) := \mu(\Pi_F^{-1}(A)), $$
where $\Pi_F$ denotes the orthogonal projection onto $F$ and $\Pi_F^{-1}(A) := \{x \in \R^n ; \Pi_F(x) \in A\}$. 

Notice that $\Pi_F^{-1}(\e^t A) = \e^t (A \times F^{\perp})$, where $F^{\perp}$ denotes the orthogonal complement of $F$. Hence if $\mu$ satisfies the (B)-property, then $\Pi_F \mu$ satisfies the (B)-property. \\

{\bf Step 2:} \underline{Approximation of log-concave measures} \\
Let us show that for every compactly supported log-concave measure $\mu$ in $\R^n$ there exists a sequence $(K_p)_{p \in \N^*}$ of convex subsets of $\R^{n+p}$ such that $\lim_{p \to + \infty} \Pi_{\R^n} \mu_{K_p} = \mu$ in the sense that the density function of $\mu$ is the pointwise limit of the density functions of $(\mu_{K_p})_{p \in \N^*}$, where $\mu_{K_p}$ denotes the uniform measure on $K_p$ (up to a constant).

Let $\mu$ be a compactly supported log-concave measure in $\R^n$ with density function $f = \e^{-V}$, where $V : \R^n \to \R \cup \{+ \infty\}$ is a convex function. Notice that for every $x \in \R^n$, $\e^{-V(x)} = \lim_{p \to +\infty} (1 - \frac{V(x)}{p})_+^p$, where for every $a \in \R$, $a_+ = \max(a,0)$. Let us define for every $p \in \N^*$,
$$ K_p = \{(x,y) \in \R^n \times \R^p ; |y| \leq \left( 1 - \frac{V(x)}{p} \right)_+\}. $$
One has for every $x \in \R^n$,
\begin{eqnarray*}
\left( 1 - \frac{V(x)}{p} \right)_+^p & = & \int_0^{\left( 1 - \frac{V(x)}{p} \right)_+} pr^{p-1} \, \de r \\ & = & p \int_0^{+ \infty} 1_{[0, \left( 1 - \frac{V(x)}{p} \right)_+]}(r) r^{p-1} \, \de r \\ & = & \frac{1}{v_p} \int_{\R^p} 1_{K_p}(x,y) \, \de y,
\end{eqnarray*}
the last inequality follows from an integration in polar coordinates, where $v_p$ denotes the volume of the Euclidean closed unit ball in $\R^p$. By denoting $\mu_{K_p}$ the measure in $\R^{n+p}$ with density function
$$ \frac{1}{v_p} 1_{K_p}(x,y), \quad (x,y) \in \R^n \times \R^p, $$
it follows that for every $p \in \N^*$, the measure $\Pi_{\R^n} \mu_{K_p}$ has density function
$$ \left( 1 - \frac{V(x)}{p} \right)_+^p, \quad x \in \R^n. $$
We conclude that $\lim_{p \to + \infty} \Pi_{\R^n} \mu_{K_p} = \mu$. \\

{\bf Step 3:} \underline{Conclusion} \\
Let $n \in \N^*$ and let $\mu$ be a symmetric log-concave measure in $\R^n$. By approximation, one can assume that $\mu$ is compactly supported. Since $\mu$ is symmetric, the sequence $(K_p)_{p \in \N^*}$ defined in Step~2 is a sequence of symmetric convex subsets of $\R^{n+p}$. If we assume that the (B)-property holds for all uniform measures in $\R^m$, for every $m \in \N^*$, then for every $p \in \N^*$, $\mu_{K_p}$ satisfies the (B)-property. It follows from Step~1 that for every $p \in \N^*$, $\Pi_{\R^n} \mu_{K_p}$ satisfies the (B)-property. Since $\lim_{p \to + \infty} \Pi_{\R^n} \mu_{K_p} = \mu$ ($c.f.$ Step~2) and since a pointwise limit of log-concave functions is log-concave, we conclude that $\mu$ satisfies the (B)-property.
\end{proof}

Similarly, let us prove that the functional form of the (B)-conjecture (Conjecture~\ref{function}) is equivalent to the classical (B)-conjecture (Conjecture~\ref{(B)}).

\begin{prop}\label{equiv}

One has equivalence between the following properties:
\begin{enumerate}
\item For every $n \in \N^*$, for every symmetric log-concave measure $\mu$ in $\R^n$ and for every symmetric convex subset $A$ of $\R^n$, the function $t \mapsto \mu(\e^t A)$ is log-concave on $\R$.
\item For every $n \in \N^*$, for every even log-concave functions $f,g : \R^n \to \R_+$, the function $t \mapsto \int_{\R^n} f(\e^{-t}x) g(x) \, \de x $ is log-concave on $\R$.
\end{enumerate}

\end{prop}

\begin{proof}
2. $\Longrightarrow$ 1. This is clear by taking $f=1_A$, the indicator function of a symmetric convex set $A$, and by taking $g$ to be the density function of a log-concave measure $\mu$. \\
1. $\Longrightarrow$ 2. Let $f,g : \R^n \to \R_+$ be even log-concave functions. By approximation, one may assume that $f$ and $g$ are compactly supported. Let us write $g = \e^{-V}$, where $V : \R^n \to \R \cup \{+ \infty\}$ is an even convex function. One has
$$ G(t): = \int_{\R^n} f(\e^{-t}x) \e^{-V(x)} \, \de x = \lim_{p \to + \infty} \int_{\R^n} f(\e^{-t}x) \left( 1 - \frac{V(x)}{p} \right)_+^p \, \de x, $$
where for every $a \in \R$, $a_+ = \max(a,0)$. Let us denote for $t \in \R$,
$$ G_p(t) = \int_{\R^n} f(\e^{-t}x) \left( 1 - \frac{V(x)}{p} \right)_+^p \, \de x. $$
We have seen in the proof of Proposition~\ref{reduction} that
$$ \left( 1 - \frac{V(x)}{p} \right)_+^p = \frac{1}{v_p} \int_{\R^p} 1_{K_p}(x,y) \, \de y, $$
where $K_p := \{(x,y) \in \R^n \times \R^p ; |y| \leq \left(  1 - \frac{V(x)}{p} \right)_+\}$ and where $v_p$ denotes the volume of the Euclidean closed unit ball in $\R^p$. Hence,
$$ G_p(t) = \frac{1}{v_p} \int_{K_p} f(\e^{-t} x) 1_{\R^p}(y) \, \de x \, \de y. $$
Notice that $K_p$ is a symmetric convex subset of $\R^{n+p}$. The change of variable $\widetilde{x} = \e^{-t}x$ and $\widetilde{y} = \e^{-t}y$ leads to
$$ G_p(t) = \frac{\e^{t(n+p)}}{v_p} \mu_p(\e^{-t} K_p), $$
where $\mu_p$ is the measure with density function
$$ h(x,y) = f(x)1_{\R^p}(y), \quad (x,y) \in \R^n \times \R^p. $$
Since a pointwise limit of log-concave functions is log-concave, we conclude that the function $G$ is log-concave on $\R$ as the pointwise limit of the log-concave functions $G_p$, $p \in \N^*$.
\end{proof}

Recall that the (B)-conjecture holds true for the Gaussian measure and for the unconditional case (see~\cite{CFM}). It follows from Proposition~\ref{equiv} that Conjecture~\ref{function} holds true if one function is the density function of the Gaussian measure or if both functions are unconditional.

\end{document}